\UseRawInputEncoding
\documentclass[12pt]{krantz}
\usepackage{fixltx2e,fix-cm}
\usepackage{amssymb}
\usepackage{amsmath}
\usepackage{graphicx}
\usepackage{subfigure}
\usepackage{makeidx}
\usepackage{multicol}
\newcommand{\pl}{\parallel}
\newcommand{\openr}{\hbox{${\rm I\kern-.2em R}$}}

\begin{document}

\chapter{Why Machine Learning Cannot Ignore Maximum Likelihood Estimation}

\vspace{-25pt}
Mark J. van der Laan\\
University of California, Berkeley\\
\texttt{laan@berkeley.edu}\\

\noindent Sherri Rose\\
Stanford University\\
\texttt{sherrirose@stanford.edu}
\vspace{-10pt}

\section{Introduction}
\begin{shortbox}
\Boxhead{The Avalanche of Machine Learning}
\noindent The growth of machine learning as a field has been accelerating with increasing interest and publications across fields, including statistics, but predominantly in computer science. How can we parse this vast literature for developments that exemplify the necessary rigor? How many of these manuscripts incorporate foundational theory to allow for statistical inference? Which advances have the greatest potential for impact in practice? One could posit many answers to these queries. Here, we assert that one essential idea is for machine learning to integrate maximum likelihood for  estimation of functional parameters, such as prediction functions and conditional densities.
\end{shortbox}

 The statistics literature proposed the familiar maximum likelihood estimators (MLEs) for parametric models and established that these estimators are asymptotically linear  such that their $\sqrt{n}$-scaled and mean-zero-centered version is asymptotically normal \cite{stigler2007epic}.  
This allowed the construction of confidence intervals and formal hypothesis testing. However, due to the restrictive form of these models, these parametric model-based MLEs target a projection of the true density on the parametric model that is hard to interpret. In response to this concern with misspecified parametric models, a rich statistical literature has studied so-called sieve-based or sieve MLEs involving specifying a sequence of parametric models that grow towards the desired true statistical model \cite{newey1997convergence,shen1997methods}. Such a sieve will rely on a tuning parameter that will then need to be selected with some method such as cross-validation. 

Simultaneously, and to a large degree independently from the statistics literature, the computer science literature developed a rich set of tools for constructing data-adaptive estimators of functional parameters, such as a density of the data, although this literature mostly focused on learning prediction functions. 
This has resulted in a wealth of machine learning algorithms using a variety of strategies to learn the target function from the data. In addition, a more recent literature in statistics established  super learning based on cross-validation as a general approach to optimally combine a library of such candidate machine learning algorithms into an ensemble that performs asymptotically as well as an oracle ensemble under specified conditions \cite{polley2011super,van2007super}.
Although the latter framework is optimal for fast learning of the target function, it lacks formal statistical inference for smooth features of the target function and for the target function itself.

\begin{shadebox}
In this chapter, we will argue that, in order to preserve statistical inference, we should preference sieve MLEs as much as possible. If the target function is not a data density, one can use, more generally, minimum loss estimation with an appropriate loss function. For example, if the goal is to learn the conditional mean of an outcome, one could use the squared-error loss and use minimum least squares estimation. Therefore, in this chapter, the abbreviation MLE also represents the more general minimum loss estimator. 
\end{shadebox}

We will demonstrate the power of sieve MLEs with a particular sieve MLE that also relies on the least absolute shrinkage and selection operator (LASSO) \cite{tibshirani1996regression} and is termed the highly adaptive LASSO minimum loss estimator (HAL-MLE) \cite{benkeser2016highly,van2017generally}. We will  argue that the HAL-MLE is a particularly powerful sieve MLE, generally theoretically superior to other types of sieve MLEs. 
Moreover, for obtaining statistical inference for estimands that aim to equal or approximate a causal quantity defined in an underlying causal model, we discuss the combination of HAL-MLE with targeted maximum likelihood estimators (TMLEs) \cite{van2011targeted,van2018targeted,van2006targeted} in HAL-TMLEs as well as using undersmoothed HAL-MLE by itself as a powerful statistical approach \cite{ju2020robust,van2017generally,van2019efficient,van2018targeted}. HALs can be used for the initial estimator in TMLEs as well as for nuisance functions representing the treatment and missingness/censoring mechanism. 
We will also clarify why many machine learning algorithms are not suited for statistical inference by having deviated from sieve MLEs. This chapter is a compact summary of recent and ongoing research; further background and references can be found in the literature cited. Sections~\ref{sec:nonpMLE}-\ref{sec:statprop} and \ref{sec:TMLE} are more technical in nature, introducing core definitions and properties. Some readers may be interested in jumping to the material in Sections~\ref{sec:contrasting}, \ref{sec:implement}, and \ref{sec:apps} for technical but narrative discussions on contrasting HAL-MLE with other estimators as well as implementation of HAL-MLE and HAL-TMLE.

\section{Nonparametric MLE and Sieve MLE}
\label{sec:nonpMLE}
Before introducing the HAL-MLE, we describe the technical foundations of nonparametric and sieve MLEs. 

\subsection{Nonparametric MLE}
We have observational unit on which we observe a random vector of measurements $O$ with probability distribution $P_0$ (e.g., $O=(W,A,Y)$ with covariates $W$, intervention or treatment $A$, and outcome $Y$). Consider the case that we observe $n$ i.i.d. copies $O_1,\ldots,O_n$ from the common probability measure $P_0$ and suppose that we know that $P_0$ is an element of a set ${\cal M}$ of possible probability distributions, which is called the statistical model. Here, the underscore $0$ indicates the true unknown probability measure. 
In addition, suppose that our target functional parameter is defined by a mapping $Q:{\cal M}\rightarrow D[0,1]^d$ from the probability measures into the space of $d$-variate real-valued cadlag functions \cite{neuhaus1971weak}. 

Let ${\cal Q}=\{Q(P): P\in {\cal M}\}$ be the parameter space for this functional parameter. Moreover, we assume that the target functional estimand $Q(P_0)$ can be characterized as the minimizer: $Q_0=\arg\min_{Q\in {\cal Q}}P_0L(Q)$, for some loss function 
$L(Q,O)$ also denoted with $L(Q)$ when viewed as a function of $O$.  A first natural attempt for estimation of $Q_0$ would be to define the MLE over the entire parameter space:
\[
\hat{Q}(P_n)=\arg\min_{Q\in {\cal Q}}P_n L(Q),\]
where we use empirical process notation $P_nf\equiv 1/n\sum_{i=1}^n f(O_i)$.

\begin{shortbox}
\noindent However, for realistic statistical models ${\cal M}$, the parameter space is generally infinite dimensional and too flexible such that the minimizer of the empirical risk will overfit to the data and results in an ill-defined $\hat{Q}(P_n)$ or a statistically poor  estimator. 
\end{shortbox}

\subsection{Sieve MLE}
A sieve MLE recognizes this need to regularize the nonparametric MLE, and proposes a sequence ${\cal Q}_{\lambda}\subset{\cal Q}$ of subspaces  indexed by some Euclidean-valued tuning parameter $\lambda \in {\cal I}\subset \openr^k$ that measures the complexity of the space ${\cal Q}_{\lambda}$ \cite{newey1997convergence,shen1997methods}. In addition, the sieve is able to approximate the target function such that for a sequence of values $\lambda_j$, ${\cal Q}_{\lambda_j}$ approximates  the entire parameter space ${\cal Q}$ as $j$ approximates infinity. 
These spaces are restricted in a manner so that an MLE: $Q_{n,\lambda}=\arg\min_{Q\in {\cal Q}(\lambda)}P_n L(Q)$, is well defined for all $\lambda$, or, at least, for a range ${\cal I}_n$ of $\lambda$ values that approximates the complete index set ${\cal I}$ as sample size $n$ grows. 
We can then define a cross-validation selector $\lambda_{n,cv}=\arg\min_{\lambda}1/V\sum_v P_{n,v}^1 L(\hat{Q}_{\lambda}(P_{n,v}))$ involving a $V$-fold sample splitting of the sample of $n$ observations where the $v$-th sample split yields an empirical distribution $P_{n,v}^1$ and $P_{n,v}$ of the validation and training sample, respectively. 
Uniform bounds on this loss function implies that this cross-validation selector performs as well as the oracle selector: $\lambda_{0,n}=\arg\min_{\lambda}1/V\sum_v P_0 L(\hat{Q}_{\lambda}(P_{n,v}))$ \cite{van2003unified,van2006cross,van2006oracle}.

\subsection{Score Equations Solved by Sieve MLE}
\begin{shortbox}
\noindent A key property of a sieve MLE is that it solves score equations. That is, one can construct a family of one dimensional paths $\{Q_{n,\lambda,\delta}^h:\delta\}\subset {\cal Q}$ through $Q_{n,\lambda}$ at $\delta=\delta_0\equiv 0$ indexed by a direction $h\in {\cal H}$, and we have:
\[
0=\frac{d}{d\delta_0}P_n L(Q_{n,\lambda,\delta_0}^h)\mbox{ for all paths $h\in {\cal H}$.}\]
\end{shortbox}
\noindent We can define a path-specific score: $A_{Q_{n,\lambda}}(h)\equiv \frac{d}{d\delta_0}L(Q_{n,\lambda,\delta_0})$, for all $h$. This mapping $A_{Q_{n,\lambda}}(h)$ will generally be a linear operator in $h$ and can therefore be extended to a linear operator on a Hilbert space generated or spanned by these directions ${\cal H}$.
One will then also have that $P_nA_{Q_{n,\lambda}}(h)=0$ for any  $h$ in the linear span of ${\cal H}$, thereby obtaining that $Q_{n,\lambda}$ solves a space of scores. For scores $S_{Q_{n,\lambda}}$ that can be approximated by scores $A_{Q_{n,\lambda}}(h)$ for certain $h$, this might provide the basis for showing that   $P_n S_{Q_{n,\lambda}}=o_P(n^{-1/2})$ is solved up to the desired approximation. 
In particular, one can apply this class of score equations at $\lambda=\lambda_{n,cv}$ to obtain the score equations solved by the cross-validated sieve MLE $Q_{n,\lambda_{n,cv}}$.

\section{Special Sieve MLE: HAL-MLE}
\label{sec:HALMLE}

We now introduce a special sieve MLE, the HAL-MLE. The HAL-MLE is generally theoretically superior to other types of sieve MLEs, with its statistical properties presented in the subsequent section.
\subsection{Definition of HAL-MLE}
\begin{shadebox}
\noindent A particular sieve MLE is given by the HAL-MLE. It starts with a representation $Q(x)=Q(0)+\sum_{s\subset\{1,\ldots,d\}}\int_{(0,x_s]}Q(du_s,0_{-s})$, where
$x_s=(x(j): j\in s)$, $0_{-s}=(0(j): j\not \in s)$ and $Q(du_s,0_{-s})$ is a generalized difference for an infinitesimal small cube $(u_s,u_s+du_s]\subset\openr^{\mid s\mid}$ in the $s$-dimensional Euclidean space with respect to the measure generated by the $s$-specific section $Q_s(u_s)\equiv Q(u_s,0_{-s})$ that sets the coordinates outside $s$ equal to zero \cite{benkeser2016highly,gill1995inefficient,van2017generally}.
\end{shadebox}
\noindent This representation shows that any cadlag function can be represented as a linear combination of indicator basis functions $x\rightarrow I_{u_s\leq x_s}$ indexed by knot point $u_s$, across all subsets $s$ of $\{1,\ldots,d\}$ with coefficients given by $Q(du_s,0_{-s})$. In particular, if $Q$ is such that all its sections are like discrete cumulative distribution functions, then this representation of $Q$ is literally a finite linear combination of these indicator basis functions. We note that each basis function $I_{u_s\leq x_s}$ is a tensor product $\prod_{j\in s}I(u_j\leq x_j)$ of univariate indictor basis functions $I(u_j\leq x_j)$, functions that jump from $0$ to $1$ at knot point $u_j$. Moreover, the $L_1$-norm of these coefficients in this representation is the so-called sectional variation norm:
\[
\pl Q\pl_v^*=\mid Q(0)\mid+\sum_{s\subset\{1,\ldots,d\}}\mid Q(du_s,0_{-s})\mid.\]

This sectional variation norm of $Q$ represents a measure of the complexity of $Q$. 
Functions that can be represented as $Q(x)=\sum_{j=1}^d Q_j(x_j)$ would correspond with 
$Q(du_s,0_{-s})=0$ for all subsets $s$ of size larger than or equal to $2$. Such functions would have a dramatically smaller sectional variation norm than a function that requires high-level interactions. If the dimension $d=1$, then this sectional variation norm is just the typical variation norm of a univariate function obtained by summing the absolute values of its changes across small intervals. For dimensional $d=2$, this is defined the same way, but now one sums so-called generalized differences $Q(dx_1,dx_2)=Q(x_1+dx_1,x_2+dx_2)-Q(x_1+dx_1,x_2)-Q(x_1,x_2+dx_2)+Q(x_1,x_2)$, representing the measure $Q$ assigns to the small 2-dimensional cube $(x,x+dx]$. Similarly, one defines generalized differences for general $d$-dimensional function as the measure it would give to a $d$-dimensional cube $(x,x+dx]$, treating the function as a $d$-variate cumulative distribution function. 

This measure of complexity now defines  the sieve ${\cal Q}_{\lambda}=\{Q\in {\cal Q}: \pl Q\pl_v^*\leq \lambda\}$ as all cadlag functions with a sectional variation norm smaller than $\lambda$. The sieve-based MLE is then defined as: 
$Q_{n,\lambda}=\arg\min_{Q:\pl Q\pl_v^*\leq \lambda}P_nL(Q)$. Using a discrete support $\{u_{s,j}:j\}$  for each $s$-specific section $Q(u_s,0_{0_{-s}})$, we can represent $Q=Q_{\beta}=\sum_{s,j}\beta(s,j)\phi_{s,j}$ with $
\phi_{s,j}(x)=I(x_s\leq u_{s,j})$, where $u_{s,j}$ represent the knot points. Thus, we  obtain a standard LASSO estimator:
\[
\beta_{n,\lambda}=\arg\min_{\beta,\pl \beta\pl_1\leq \lambda}P_n L(Q_{\beta}),\]
and resulting $Q_{n,\lambda}=Q_{\beta_{n,\lambda}}$.
As above, $\lambda$ is then selected with the cross-validation selector $\lambda_{n,cv}$. 
This can generally be implemented with LASSO software implementations such as \texttt{glmnet} in \texttt{R} \cite{friedman2021package}. Additionally, \texttt{HAL9001} provides  implementations for linear, logistic, Cox, and Poisson regression, which also provides HAL estimators of the conditional hazards and intensities \cite{hejazi2020hal9001}. 

\subsection{Score Equations Solved by HAL-MLE}
In order to define a set of score equations solved by this HAL-MLE we can use as paths
$\beta_{\delta}^h(s,j)=(1+\delta h(s,j))\beta(s,j)$ for a vector $h$. In order to keep the $L_1$-norm constant along this path, we then enforce the constraint $r(h,\beta)=0$, where
\[
r(h,\beta)\equiv h(0)\mid \beta(0)\mid+\sum_{s,j} h(s,j)\mid \beta(s,j)\mid.\]
It can be verified that now $\pl \beta_{\delta}^h\pl_1=\pl \beta\pl_1$ for  $\delta$ in a neighborhood of $0$.
Therefore, we know that the HAL-MLE solves the score equation for each of these paths: 
\[
P_n A_{\beta_{n,\lambda}}(h)=0\mbox{ for all $h$ with $r(h,\beta_{n,\lambda})=0$.}\]
This shows that the HAL-MLE  solves a class of score equations. Moreover, this result can be used to prove that 
\[P_n \frac{d}{d\beta_{n,\lambda}(j)}L(Q_{\beta_{n,\lambda}})=O_P(n^{-2/3}),\] for all $j$ for which $\beta_{n,\lambda}(j)\not =0$. That is, the $L_1$-norm constrained HAL-MLE also solves the unconstrained scores solved by the MLE over the finite linear model $\{Q_{\beta}:\beta(j)=0$ if $\beta_{n,\lambda}(j)=0\}$ implied by the fit $\beta_{n,\lambda}$. 
By selecting the $L_1$-norm $\lambda$ to be larger, this set of score equations approximates any score, thereby establishing that the HAL-MLE, if slightly undersmoothed, will solve  score equations up to $O_P(n^{-2/3})$ uniformly over all scores. For further details we refer prior work \cite{van2019efficient}.

\begin{shortbox}
\noindent This capability of the HAL-MLE to solve all score equations, even uniformly over a class that will contain any desired efficient influence curve of a target parameter, provides the fundamental basis for establishing its remarkably strong asymptotic  statistical performance in estimation of smooth features of $Q_0$ as well as of $Q_0$ itself. 
\end{shortbox}

\section{Statistical Properties of the HAL-MLE}
\label{sec:statprop}

Having introduced the special sieve MLE, the HAL-MLE $Q_{n,\lambda_n}$, we now further enumerate its statistical properties.

\subsection{Rate of Convergence}
Firstly, 
$d_0(Q_{n,\lambda_n},Q_0)=O_P(n^{-2/3}(\log n)^d)$, where $d_0(Q,Q_0)=P_0L(Q)-P_0L(Q_0)$ is the loss-based dissimilarity. This generally establishes that some $L^2$-norm of $Q_{n,\lambda_n}-Q_0$ will be $O_P(n^{-1/3}(\log n)^{d/2})$ \cite{bibaut2019fast,van2017generally}. 

\subsection{Asymptotic Efficiency for Smooth Target Features of Target Function}
Let $\Psi(Q_0)$ be a pathwise differentiable target feature of $Q_0$ and let $D^*_{Q(P),G(P)}$ be its canonical gradient at $P$, possibly also depending on a nuisance parameter $G(P)$.  
Let $R_{P_0}(Q,Q_0)=\Psi(Q)-\Psi(Q_0)+P_0 D^*_{Q,G_0}$ be its exact second-order remainder.  
By selecting $\lambda_n$ to be large enough so that $P_n D^*_{Q_{n,\lambda_n},G_0}=o_P(n^{-1/2})$, we obtain that $\Psi(Q_{n,\lambda_n})$ is asymptotically linear with influence curve $D^*_{Q_0,G_0}$. That is, $\Psi(Q_n)$ is an asymptotically efficient estimator of $\Psi(Q_0)$. The only conditions are that $R_{P_0}(Q_{n,\lambda_n},Q_0)=o_P(n^{-1/2})$, where one can use that $d_0(Q_{n,\lambda_n},Q_0)=O_P(n^{-2/3}(\log n)^d)$. Generally speaking, the latter will imply that the exact remainder has the same rate of convergence. In fact, in double robust estimation problems, this exact remainder equals zero (because it is at $G_0$).

\subsection{Global Asymptotic Efficiency}
Moreover, now consider  a very large class of pathwise differentiable target features $\Psi_t(Q_0)$ indexed by some $t\in [0,1]^m$, where $D^*_{t,Q,G}$ is the canonical gradient and $R_{P_0,t}(Q,Q_0)$ is the exact remainder.  Then, by undersmoothing enough so that
$\sup_t \mid P_n D^*_{t, Q_{n,\lambda_n},G_0}\mid=o_P(n^{-1/2})$, that is, we use global undersmoothing, we obtain that
$\Psi_t(Q_n)$ is an asymptotically linear estimator of $\Psi_t(Q_0)$ having influence curve $D^*_{t,Q_0,G_0}$, with a remainder  $R_n(t)$ that satisfies $\sup_t\mid R_n(t)\mid =o_P(n^{-1/2})$, so that $n^{1/2}(\Psi_t(Q_n)-\Psi_t(Q_0))$ converges weakly to a Gaussian process as a random function in function space. That is, $(\Psi_t(Q_n):t\in [0,1]^m)$ is an asymptotically efficient estimator of $(\Psi_t(Q_0):t\in [0,1]^m)$ in a sup-norm sense. 
In particular, it allows construction of simultaneous confidence bands. We note that this also implies  that $Q_{n,h}(x)\equiv \int Q_n(u)K_h(u-x)du$, a kernel smoother of $Q_n$ at $x$ is an asymptotically efficient estimator of $Q_{0,h}(x)\equiv \int Q_0(u)K_h(u-x)du$, a kernel smoother of the true target function $Q_0$. In fact, this holds uniformly in all $x$, and for a fixed $h$ we have weak convergence of $n^{1/2}(Q_{n,h}-Q_{0,h})$ to a Gaussian process.  That is, the undersmoothed HAL-MLE is an efficient estimator of the kernel smoothed functional of $Q_0$, for any $h$. 

\begin{extract}
\textit{This may make one wonder if $Q_n(x)$ itself is not asymptotically normally distributed as well? While not a currently solved problem,  we conjecture that, indeed, $(Q_n(x)-Q_0(x))/\sigma_n(x)$ converges weakly to a normal distribution, where the rate of convergence may be as good as  $n^{-1/3}(\log n)^{d/2}$. 
If this results holds, then the HAL-MLE also allows formal statistical inference for the function itself! }
\end{extract}

\subsection{Nonparametric Bootstrap for Inference}
It has been established  that the nonparametric bootstrap consistently estimates the multivariate normal limit distribution of the plug-in HAL-MLE for a vector of target features  \cite{cai2020nonparametric}. 
In fact, it will consistently estimate the limit of a Gaussian process of the plug-in HAL-MLE for an infinite class of target features as discussed above. 
The approach is computationally feasible and still correct by only bootstrapping the LASSO-based fit of the model that was selected by the HAL-MLE; the bootstrap is only refitting a high-dimensional linear regression model with at most $n$ basis functions (and generally much smaller than $n$). 
This allows one to obtain inference that also attains second-order behavior. In particular,  obtaining the bootstrap distribution for each $L_1$-norm larger than or equal to the cross-validation selector of the $L_1$-norm and selecting the $L_1$-norm at which the confidence intervals plateau provides a robust {\em finite sample} inference procedure \cite{cai2020nonparametric}.

It will also be of interest to understand the behavior of the nonparametric bootstrap in estimating the distribution of the HAL estimator of the target function itself (say, at a point). Before formally addressing this, we need to first establish that the HAL-MLE is asymptotically normal as conjectured in a previous subsection.

\subsection{Higher-Order Optimal for Smooth Target Features}
In recent work, higher-order TMLEs were proposed that target extra score equations beyond the first-order canonical gradient of the target feature, which are selected so that  the exact second-order remainder for the plug-in TMLE will become a higher-order difference (say, third-order for the second order TMLE) \cite{van2021higher}.  Because an undersmoothed HAL-MLE also solves these extra score equations, it follows that the plug-in HAL-MLE, when using some undersmoothing, is not just asymptotically efficient but also has a reduced exact remainder analogue to the higher-order TMLE. (More on TMLEs in later sections.)

The key lesson is that the exact remainder in an expansion of a plug-in estimator can be represented as an expectation of a score: $R_0((Q,G),(Q_0,G_0))=P_0 A_Q(h_{P,P_0})$, for some direction $h_{P,P_0}$. For example, in the nonparametric estimation of $E(Y_1)=E[E(Y\mid A=1,W)]$, where $\bar{Q}=E(Y|A=1,W)$ and $\bar{G}(W)=E(A\mid W)$, we have
$R_0((Q,G),(Q_0,G_0))=P_0 (\bar{Q}-\bar{Q}_0)(\bar{G}-\bar{G}_0)/\bar{G}$, which can be written as $-E_0 I(A=1)/\bar{G}_0(W)(\bar{G}-\bar{G}_0)/\bar{G}(Y-\bar{Q})$, which is an expectation of a score of the form $h_{G,G_0}(A,W)(Y-\bar{Q})$. The HAL-MLE $\bar{Q}_n$ approximately solves $P_n h_{G_n,G_0}(Y-\bar{Q}_n)$ for some approximation $G_n$ of $G_0$, thereby  reducing the remainder to
$(P_n-P_0)h_{G_n,G_0}(Y-\bar{Q}_n)$. The latter is a term that is $o_P(n^{-1/2})$---only relying on the consistency of $G_n$---not even a rate is required.  Higher-order TMLE directly targets these scores to map the remainder in higher and higher order differences, but the undersmoothed HAL-MLE is implicitly doing the same (although not solving the score equations exactly).

\begin{shortbox}
\noindent This demonstrates that solving score equations has enormous implications for first- and higher-order behavior of a plug-in estimator. Typical machine learning algorithms are generally not tailored to solve score equations, and, thereby, will not be able to achieve such statistical performance for their plug-in estimator. In fact, most machine learning algorithms are not grounded in any asymptotic limit distribution theory.
\end{shortbox}

\section{Contrasting HAL-MLE with Other Estimators}
\label{sec:contrasting}
One might wonder: what is particularly special about this  sieve in the HAL-MLE relative to other sieve MLEs, such as those using Fourier basis or polynomial basis, wavelets, or other sequences of parametric models. Also, why might we prefer it over other general machine learning algorithms?

\subsection{HAL-MLE vs.\ Other Sieve MLE}

\begin{shadebox}
\noindent The simple answer is that these sieve MLEs generally do not have the (essentially) dimension-free/smoothness-free rate of convergence $n^{-1/3}(\log n)^{d/2}$, but instead their rates of convergence heavily rely on assumed smoothness. HAL-MLE uses a global complexity property, the sectional variation norm, rather than relying on local smoothness. The global bound on the sectional variation norm provides a class of cadlag functions  ${\cal F}$ that has a remarkably nice covering number $\log N(\epsilon,{\cal F},L^2)\lesssim 1/\epsilon (\log (1/\epsilon))^{2d-1}$, hardly depending on the dimension $d$. Due to this covering number, the HAL-MLE has this powerful rate of convergence---only assuming that the true target function is a cadlag function with finite sectional variation norm. 
\end{shadebox}

A related advantage of the special HAL sieve   is that the union of all indicator basis functions is a small Donsker class, even though it is able to span any cadlag function. Most sets of basis functions include ``high frequency'' type basis functions that have a variation norm approximating infinity. As a consequence, these basis functions do not form a nice Donsker class. In particular, this implies that the HAL-MLE does not overfit, as long as the sectional variation norm is controlled.  The fact that HAL-MLE itself is situated in this  Donsker class also means that the efficient influence curves at HAL-MLE fits will  fall in a similar size Donsker class. As a consequence, the Donsker class condition for asymptotic efficiency of plug-in MLE and TMLE is naturally satisfied when using the HAL-MLE, while other sieve-based estimators easily cause a violation of the Donsker class condition.

\begin{shortbox}
\noindent This same powerful property of the Donsker class spanned by these indicator basis functions also allows one to prove that nonparametric bootstrap works for the HAL-MLE, while, generally speaking, the nonparametric bootstrap generally fails to be consistent for machine learning algorithms. 
\end{shortbox}

Another appealing feature of the HAL sieve is that it is only indexed by the $L_1$-norm, while many  sieve MLEs rely on a precise specification of the sequence of parametric models that grow in dimension. It should be expected that the choice of this sequence can have a real  impact in practice. HAL-MLE does not rely on an ordering of basis functions, but rather it just relies on  a complexity measure. For each value of the complexity measure, it includes all basis functions and represents  an infinite dimensional class of functions rich enough to approximate any function with complexity satisfying this bound. A typical sieve MLE can only approximate the true target function, while the HAL-MLE includes the true target function when the sectional variation norm bound exceeds the sectional variation norm of the true target function. 

As mentioned, the HAL-MLE solves the regular score equations from the data-adaptively selected HAL-model at rate $O_P(n^{-2/3})$. As a consequence, the HAL-MLE is able to uniformly solve the class of all score equations---only restricted by some sectional variation norm bound, where this bound can go to infinity as sample size increases. This strong capability for solving score equations appears to be unique for HAL-MLE relative to other sieve-based estimators. It may  be mostly due to actually being an MLE over an infinite function class (for a particular  variation norm bound). We also note that a parametric model-based sieve MLE would be forced to select a small dimension to avoid overfitting. However, the HAL-MLE adaptively selects such a model among all possible basis  functions, and the dimension of this data-adaptively selected model will generally be larger. The latter is due to the HAL-MLE only being an $L_1$-regularized MLE for that adaptively selected parametric model, and thereby does not overfit the score equations, while the typical sieve MLE would represent a full MLE for the selected model.  

\begin{shortbox}
\noindent By solving many more score equations \textit{approximately} HAL-MLE can span a much larger space of scores than a sieve MLE that solves many fewer score equations perfectly. 
\end{shortbox}

\subsection{HAL-MLE vs.\ General Machine Learning Algorithms}
Sieve MLEs solve score equations and, thereby, are able to approximately solve a class of score equations---possibly enough to approximate  efficient influence curves of the target features of interest, and thus be asymptotically efficient for these target features. We discussed how these different sieve MLE can still differ in their performance in solving score equations, and that the HAL-MLE appears to have a unique strategy by choice of basis functions and measure of complexity. Thus, this gives it a distinct capability to solve essentially all score equations at a desired approximation error. 

\begin{shortbox}
\noindent  Many machine learning algorithms fail to be an MLE over any subspace of the parameter space. Such algorithms will have poor performance in solving score equations. As a consequence, they will not result in asymptotically linear plug-in estimators and will generally be overly biased and nonrobust. 
\end{shortbox}

\section{Considerations for Implementing HAL-MLE}
\label{sec:implement}
Having presented core definitions, properties, and comparisons, we now turn to some additional considerations for implementing HAL-MLEs.

\subsection{HAL-MLE Provides Interpretable Machine Learning}
The role of interpretable algorithms is a major consideration for many applications \cite{rudin2019stop}. Is HAL-MLE interpretable? The HAL-MLE is a finite linear combination of indicator basis functions, so-called tensor products of zero-order splines, and the number of basis functions is generally significantly smaller than $n$. HAL-MLE has also been generalized to higher-order  splines so that it is forced to have a certain smoothness. These higher-order spline HAL-MLE have the same statistical properties as reviewed above, as long as the true function satisfies the enforced level of smoothness.  Such smooth HAL-MLE can be expected to result in even sparser fits, e.g., a smooth monotone function can be fitted with a few first-order splines (piecewise linear) while it needs relatively many knot points when fitting with a piecewise constant function.  
Either way, the HAL-MLE is a closed form object that can be evaluated and is thus interpretable. Therefore, the HAL-MLE has the potential to play a key role in interpretable machine learning. 

\subsection{Estimating the Target Function Itself}
\begin{shadebox}
\noindent  HAL-MLE requires selecting the set of knot points and, thereby, the collection of spline basis functions. The largest set of knot points we have recommended (and suffices for obtaining the nonparametric rate of convergence $n^{-1/3}(\log n)^{d/2})$ is given by:  $\{X_{s,i}: i=1,\ldots,n,s\subset\{1,\ldots,d\}\}$, which corresponds (for continuous covariates) to $N=n (2^d-1)$ basis functions. 
The LASSO-based fit will then select a relatively small subset ($<<n$) of this user-supplied collection. 
\end{shadebox}

Rather than selecting this full set of basis functions, one can  incorporate model assumptions. This could include only selecting up to two-way tensor products, ranking the basis functions by their sparsity  (i.e., proportion of 1's among the $n$ observations) and selecting the top $k$, and specifying a specific additive model using standard \texttt{glm}-formula notation, such as 
$f(x_1,x_2,x_3)=f_1(x_1)+f_2(x_2,x_3)+f_3(x_3)$, and selecting the knot points accordingly. In particular, one could use some \texttt{glmnet} fit using main terms and standard interactions to decide on this type of additive model. In addition, one can specify that the coefficients of a certain set of basis functions should be nonnegative and others should be non-positive, thereby enforcing monotonicity of functions in the additive model.  Finally, one can select among zero order and more generally $k$-th order spline basis functions, thereby specifying a desired smoothness of the HAL-MLE.

\begin{shortbox}
\noindent  To reduce the computational burden of the implementation of HAL-MLE, one can randomly subset from a large set of basis functions or subset in a deterministic manner. 
\end{shortbox}

For example, for a continuous covariate $X_j$, rather than using as knot points $X_{j,i}$, $i=1,\ldots,n$, one selects only 5 knot points---the 5 quantiles of the $n$ observations $X_{j,i}$. In this manner, the number of basis functions will not grow linearly in $n$, it is now growing by a fixed factor $5$.  Similarly, the above restriction on the degree of the tensor products to only $2$ reduces the $2^d$-factor to $d^2$. Therefore, such choices reduces the total number of basis functions from $n (2^d-1)$ to  around $5*d^2$ basis functions. 

In certain cases, one might view some of these restrictions as a specification of the actual statistical model. For example, the statistical model might assume that $Q_0$ is an additive model including any one-way, two-way, and three-way function. 
However, in general, many of these model choices for the HAL-MLE will be  hard to defend based on prior knowledge, although they might result in a statistically improved estimator relative to using the most nonparametric HAL-MLE. Therefore, we recommend building a super learner whose library contains a variety of such HAL-MLE specifications. In addition, by ranking these HAL-MLE fits by their complexity, one could implement a cross-validation scheme that implements and evaluates estimators from least complex to increasingly complex and stops when the cross-validated risk of the next HAL-MLE drops below the performance of the less complex previous HAL-MLE. In this manner, the resulting super learner avoids having to implement the highly computer intensive HAL-MLEs, except when they are really needed. Because the discrete super learner is asymptotically equivalent with the oracle selector, the resulting discrete super learner will perform as well as the best possible HAL-MLE, thereby also inheriting the rate of convergence of the most nonparametric HAL-MLE.

\subsection{Using the HAL-MLE for Smooth Feature Inference}
By specifying the class of target features for which we want the HAL-MLE to be efficient, we could undersmooth the HAL-MLE until all the efficient score equations (one for each target feature) are solved at the desired level. In this case, the above discrete super learner is not making the right trade-off. One might still be able to solve all the score equations under a restricted model choice, such as only including up to three-way tensor product basis functions. Thus, one might still rank various HAL-MLE specifications from least complex to maximal complex, but  would now keep increasing complexity of the HAL-MLE according to this ranking until the desired class of score equations are solved at the desired level (i.e., $o_P(n^{-1/2})$ uniformly in all target features).

\section{HAL and TMLE}
\label{sec:TMLE}

Given the bevy of statistical properties for the HAL-MLE, one might wonder why one would need TMLE \cite{van2011targeted,van2006targeted}? However, the finite sample statistical properties of the HAL-MLE rely on solving the desired score equations up to a small enough error, even though the HAL-MLE has the statistical capacity of solving essentially all score equations up to an asymptotic rate $O_P(n^{-2/3})$. On the other hand, if the HAL-MLE is too undersmoothed, then its performance in estimating the target function deteriorates, and that can also have a detrimental effect on  statistical  performance. Targeted undersmoothing can deal with this to some degree by precisely finding the level of undersmoothing that solves the desired target equations at the preferred level (which could be $O_P(\sigma_n/(n^{1/2}\log n))$, where $\sigma_n$ is a standard error of the efficient influence curve of the target feature). Nonetheless, in finite samples it might simply not be possible to achieve this goal. 

\subsection{Targeting a Class of Target Features}

Therefore, to achieve improved finite sample performance of the HAL-MLE for a desired target feature or class of target features, it follows to use TMLE with an initial HAL-MLE estimator \cite{van2017generally}. In this manner, the user has precise control over a set of target score equations while the HAL-MLE will provide important benefits by solving a large class of score equations and having a good rate of convergence. 

\begin{shadebox}
\noindent  Let $Q_n^0$ be a HAL-MLE, possibly a discrete super learner based on a library of different HAL-MLE specifications. Let $(\Psi_t(P_0): t)$ be our class of target features and let $(D^*_{t,Q,G}:t)$ be the corresponding class of canonical gradients. We can then construct a universal least favorable path $\{Q_{n,\epsilon}^0:\epsilon\}$ through $Q_n^0$ at $\epsilon =0$, possibly using an initial estimator $G_n$, and resulting MLE: $\epsilon_n=\arg\min_{\epsilon}P_n L(Q_{n,\epsilon}^0)$, so that
$\sup_t\mid P_n D^*_{t,Q_{n,\epsilon_n}^0,G_n}\mid =o_P(n^{-1/2})$.
Such a TMLE update step $Q_n^*=Q_{n,\epsilon_n}^0$ satisfies the identity:
\[
\Psi_t(Q_n^*)-\Psi_t(Q_0)=(P_n-P_0)D^*_{t,Q_n^*,G_n}+R_{0,t}((Q_n^*,G_n),(Q_0,G_0)).\]
\end{shadebox}

\noindent Because $d_0(Q_n^*,Q_0)=O_P(n^{-2/3}(\log n)^{d})$, and assuming a HAL-MLE $G_n$ also has $d_{01}(G_n,G_0)=O_P(n^{-2/3}(\log n)^{d})$, under a strong positivity assumption, we will have:
\[
\sup_t \mid R_{0,t}((Q_n^*,G_n),(Q_0,G_0))\mid =O_P(n^{-2/3}(\log n)^d).\]
In addition, because HAL-MLEs fall in the well understood Donsker class of cadlag functions with a uniform bound on their sectional variation norm \cite{bibaut2019fast}, it generally also follows that $\{D^*_{t,Q,G}:Q\in {\cal Q},G\in {\cal Q}\}$ is a Donsker class with the same covering number rate as this cadlag function class. Therefore, we will also have:
\[
(P_n-P_0)D^*_{t,Q_n^*,G_n}=(P_n-P_0)D^*_{t,Q_0,G_0}+O_P(n^{-2/3}(\log n)^d),\]
where the rate of the remainder would follow from using the empirical process finite sample asymptotic equicontinuity results \cite{bibaut2019fast}.
As a consequence, we have:
\[
\Psi_t(Q_n^*)-\Psi_t(Q_0)=P_n D^*_{t,Q_0,G_0}+O_P(n^{-2/3}(\log n)^d).\]
In particular, it follows that $\Psi(Q_n^*)$ is an asymptotically efficient estimator of $\Psi(Q_0)$, possibly viewed as elements in function space endowed with a supremum norm \cite{van2017generally}. 

\begin{shortbox}
\noindent Therefore, in great generality, we have that  TMLEs that use HAL as an initial estimator are asymptotically efficient for the target features targeted by the TMLE assuming only that the nuisance functions are cadlag, have finite sectional variation norm, and a strong positivity assumption. 
\end{shortbox}

\noindent  It is not required that the HAL-MLE is undersmoothed for  TMLE. The solving of the score equations is carried out by the TMLE so that the HAL-MLE can be optimized for estimating $Q_0$ and $G_0$. In particular, one can now use the discrete super learner discussed above with a library of HAL-MLEs.

\subsection{Score Equation Preserving TMLE Update}
Suppose that the HAL-MLE was undersmoothed and thereby was itself already a great estimator for these target features and also for other target features that were not targeted by the TMLE. In that case, could the TMLE  destroy the score equations solved by the HAL-MLE and  deteriorate some of the good performance of the initial HAL-MLE, such as its properties in estimation of other features or its higher-order properties in estimation of the actual target features targeted by the TMLE? That is, given  an initial estimator  has already succeeded in solving a large class of important score equations, making it potentially globally efficient and higher-order efficient for a large class of features,  we should be  worried about the TMLE not preserving these score equations during its TMLE update step solving the desired target score equations up to user supplied error. 

\begin{shortbox}
\noindent This motivates us to generalize the TMLE update step to a TMLE update that is not only solving the target score equations but also preserves the score equations already solved by the initial estimator. 
\end{shortbox}

\noindent In particular, this is a motivation for using universal least favorable paths in the definition of the TMLE update, because such paths are only maximizing the likelihood in the direction  of the target score equations, thereby not affecting any score equation orthogonal to these target score equations. However, in addition, one might do the following. 
We already specified the score equations exactly solved by the HAL-MLE above, one score for each coefficient that has a nonzero coefficient corresponding with a path that keeps the $L_1$-norm constant. 
Given this specified set of score equations and its linear span $H_n$ of scores, we could compute the projection of the first-order canonical gradient $D^*_{t,Q,G}$ onto the space $H_n$ spanned by these scores, and subtract it from $D^*_{t,Q,G}$, resulting in an orthogonalized $\tilde{D}^*_{t,Q,G}=D^*_{t,Q,G}-\Pi(D^*_{t,Q,G}\mid H_n)$.

One now defines the TMLE using the universal least favorable path based on this orthogonalized set of scores $(\tilde{D}^*_{t,Q,G}:t)$ rather than $(D^*_{t,Q,G}:t)$. In this case, the TMLE update $Q_n^*$ will still solve the score equations in $H_n$ that were solved by the initial HAL-MLE $Q_n^0$, and, in addition, it will solve the score equations $P_n\tilde{D}^*_{t,Q_n^*,G_n}=0$. As a consequence, it will also solve $P_n D^*_{t,Q_n^*,G_n}=0$. So in this way, we have used the general TMLE to obtain a TMLE that not only targets a new set of score equations but also preserves the score equations already solved by the initial estimator $Q_n^0$. 

\begin{extract}
\textit{In future work, we will study, implement, and evaluate this score equation preserving TMLE, and other variations of such score equation preserving TMLEs. The key message is that  we will have further robustified the TMLE by not only solving the target score equations and preserving the rate of convergence of the initial estimator, but also preserving the score equations solved by the initial estimator with all its important additional statistical benefits. }
\end{extract}

\section{Considerations for Implementing HAL-TMLE}
\label{sec:apps}

Causal inference relies on a general roadmap involving defining a causal model, causal quantity of interest, establishing an identification result, and thereby defining a target estimand of the data distribution \cite{pearl2009causality}. These steps then dictate the statistical estimation problem in terms of data, statistical model, and target estimand \cite{petersen2014causal,van2011targeted}. For estimation and statistical inference, the underlying causal model and its assumptions play no role, such that one now just proceeds with the roadmap of targeted learning given the statistical model and target estimand \cite{schuler2017targeted,van2011targeted}. That is, we can:  compute the canonical gradient of the target estimand, define an initial estimator, describe a TMLE update of this initial estimator, and then provide confidence intervals. In particular, the above plug-in HAL-MLE and score equation preserving TMLE using the HAL-MLE as initial estimator can now be applied. In addition, one can develop the higher-order TMLE using the HAL-MLE as initial estimator \cite{van2021higher}.  We now focus on additional considerations for implementing HAL-TMLEs.

\subsection{Double Robustness}
Many causal inference estimation problems have a double robust structure in the sense that the exact remainder $R_0((Q,G),(Q_0,G_0))$ has a cross-product structure that can be naturally bounded by a product of  $\pl Q-Q_0\pl $ and $\pl G-G_0\pl$ for certain $L^2$-norms (typically $L^2(P_0)$-norms), or, equivalently, in terms of $d_0^{1/2}(Q,Q_0)$ and $d_{01}^{1/2}(G,G_0)$. This structure implies that in randomized trials or observational studies where there is substantial knowledge of the censoring and treatment mechanism $G_0$, one can obtain excellent finite sample performance for TMLEs utilizing the resulting fast converging estimator $G_n$ of $G_0$. In such a particular appeal due to its ability to remove all bias in the initial estimator $Q_n^0$ when $G_n\approx G_0$.

\subsection{HAL for Treatment and Censoring Mechanisms}

Using HAL for $G_n$ (respecting the known model for $G_0$) has important benefits. In particular, even if $Q_n$ is inconsistent, the TMLE will remain asymptotically linear if $G_n$ is somewhat undersmoothed. In fact, if one uses a nonparametric HAL-MLE ignoring the model for $G_0$, and the full data model is nonparametric, then the TMLE will even be efficient, despite the inconsistency of $Q_n$. Similarly, the inverse probability of treatment and censoring weighted estimator using an undersmoothed HAL-MLE is asymptotically linear in these causal inference problems, and its efficiency is maximized by using an undersmoothed nonparametric HAL-MLE. On the other hand, in such a setting, if one would  use other machine learning algorithms, including a super learner, these same estimators would lose their asymptotic linearity  and not even converge at the desired $n^{-1/2}$-rate---failing to provide valid inference. 

\begin{shortbox}
\noindent Therefore, we learn that it is not only beneficial to use HAL-MLE as initial estimator of $Q_0$ due to above mentioned reasons, but it is also highly beneficial to use a HAL-MLE for $G_n$. 
\end{shortbox}

Let's consider the case where it is known that the treatment mechanism only depends on $2$ confounders. By estimating $G_0$ with a model-based HAL-MLE, perhaps undersmoothed, the above arguments show important gains. However, the above arguments also state that ignoring the model for $G_0$ will even make it more efficient. Therefore, there is an important selection problem among candidate HAL-MLEs of $G_0$ that have varied complexity, ranging from the actual model to a fully nonparametric HAL-MLE.  Cross-validation would then select the model-based HAL-MLE with high probability and thus be ignorant of the subtle bias-variance trade-offs at stake. 

More generally, selecting among candidate estimators of $G_0$ based on the log-likelihood loss can be problematic when the positivity assumption is practically violated. For example, this could result in an estimator $G_n$ that approaches zero, and, as a consequence, results in erratic TMLE updates. This is typically resolved by truncation, but one still needs to select the truncation level. Similarly, the adjustment set used by $G_n$ might need to be tailored toward the MSE of the TMLE of the target estimand using $G_n$, rather than toward the estimation of $G_0$. For example, some  baseline covariates might be highly predictive of treatment while not being real confounders (like instrumental variables), and it is well established that adjustment for such variables can easily  increase both variance and bias \cite{greenland2008invited,myers2011effects,rotnitzky2010note,schisterman2009overadjustment,schneeweiss2009high}. 
Therefore, an  important feature of causal inference problems is the targeted selection among candidate estimators of $G_0$. 

\subsection{Collaborative TMLE}
TMLE provides a natural approach for such a selection by evaluating the performance of a candidate estimator of $G_0$ with the increase of the log-likelihood during the TMLE update step. The MLE of $\epsilon$ in  the universal least favorable path through the initial estimator corresponds with fitting the target estimand itself. Specifically, it is the path along which the ratio of the squared change in target estimand by increase in log-likelihood is maximized, i.e., the slope in change of estimand per unit increase in log-likelihood is maximal. 
Therefore, this increase in log-likelihood during the TMLE update step indeed provides a targeted criterion for evaluating the performance of a candidate estimator $G_n$ of $G_0$.
After having selected an optimal choice with respect to this criterion, one can iterate the process by making the TMLE update the initial estimator and repeating this greedy selection, but now only selecting among candidate estimators that are more complex (say, higher log-likelihood) than the one selected at the previous step. In this manner, we obtain a sequence of TMLE updates $Q_{n,k}^*$ of $Q_n^0$ with corresponding increasingly complex $G_{n,k}$ estimators, $k=1,\ldots, z$. And $k$ can now be selected with cross-validation based on the loss function $L(Q)$ for $Q_0$, or a AIC/BIC type selector to avoid the double cross-validation. 
This approach can now be used, for example, to select the $L_1$-norm in a HAL-MLE, the truncation level, but also to select the maximal degree of the tensor products, or maximal sparsity of the basis functions, and so on \cite{cai2020nonparametric,ju2020robust,van2019efficient,van2021higher,van2018c}. Simulations have demonstrated that this collaborative (C-TMLE) procedure can dramatically outperform a TMLE that uses pure likelihood based estimation of $G_0$ when there are positivity issues \cite{ju2020robust}. One can also let a HAL-MLE fit of $Q_0$ impact the basis functions included in the HAL-MLE of $G_0$, and use this C-TMLE procedure to select the $L_1$-norm in the HAL-MLE of $Q_0$ and the resulting HAL-MLE of $G_0$ (as well as the truncation level). The latter type of estimator was termed the outcome adaptive HAL-TMLE \cite{ju2020robust}, which built on prior work in outcome-adaptive LASSOs \cite{shortreed2017outcome}, providing a powerful tool for variable selection in $G_n$, and was shown to have strong practical performance.

\section{Closing}

\begin{shadebox}
\noindent There is a tendency for the machine learning literature to focus on piecemeal and small extensions of the `flashy' estimator of the moment. This was recently random forests, but is currently deep learning. However, the statistical theory and empirical process literature has offered strong guidance on the development of data-adaptive estimators for features of the data distribution that provide  inference, while  fully utilizing the knowledge of a statistical model.
\end{shadebox}

 In particular, for optimal robust (higher-order)   estimation of target estimands, we need to  solve specific first- and higher-order canonical gradients of the target estimand. Also, having the functional parameter of the data distribution needed for estimation of the target estimand as a member of a function class with a  well behaving entropy integral (as implied by the covering number of the function class) provides good rates of convergence and allows for bootstrap based inference. 
HAL-MLE and HAL-TMLE appear to satisfy these key fundamental properties. 
It would be exciting for the general machine learning literature to build on these areas. 

	\bibliographystyle{plain}
	\bibliography{HAL}
\end{document}